\documentclass[12pt]{article}
\usepackage{amsmath}
\usepackage{amssymb}
\usepackage{amscd}
\usepackage[dvips]{graphics}
\usepackage{amsfonts}

\usepackage{color}
\newcommand{\refm [1]} {(\ref{#1})}

\newcounter{theorem}

\newtheorem{proposition}{Proposition}
\newtheorem{theorem}{Theorem}
\newtheorem{lemma}{Lemma}

\newcommand{\la}{\label}

\newcommand{\PP}{\mathbb{P}}

\newcommand{\nicef}{{\bf {\cal F}}}
\newcommand{\en}{\end{equation}}

\def\PP{{\Bbb P}}

\def\BR{{\Bbb R}}

\def\BE{{\Bbb E}}

\def\BZ{{\mathbb Z}}

\def\proof{\noindent{\bf Proof\ \ }}
\def\qed{\mbox{\rule{0.5em}{0.5em}}}

\def\la{\label}

\def\deli{\delta_n^{(i)}}
\def\del1{{\delta_n^{(1)}}}

\newcommand{\ignore}[1]{}

\def\BR{{\Bbb R}}

\def\BE{{\Bbb E}}
\def\BZ{{\Bbb Z}}

\def\calh{{\cal H}}

\newcommand{\be}{\begin{equation}}
\newcommand{\eu}{\end{equation}}
\newcommand{\ber}{\begin{eqnarray}}
\newcommand{\ena}{\end{eqnarray}}
\newcommand{\nin}{\noindent}
\newcommand{\non}{\nonumber}

\begin{document}

\title{A Meinardus theorem with multiple singularities}
\author{ Boris L. Granovsky\thanks{ Department of Mathematics,
Technion-Israel Institute of Technology, Haifa, 32000, Israel,
e-mail:mar18aa@techunix.technion.ac.il
 }
 \ and  Dudley
Stark\thanks{School of Mathematical Sciences, Queen Mary,
University of London, London E1 4NS, United Kingdom,
e-mail:D.S.Stark@maths.qmul.ac.uk\newline\newline 2000 {\it
5Mathematics Subject Classification}: 05A16,60F99,81T25
\newline \nin Keywords and phrases: Meinardus' theorem, Asymptotic enumeration, Dirichlet generating functions, Models of ideal gas and of quantum field theory.
} } \date{} \maketitle

\begin{abstract}
Meinardus proved a general theorem about the asymptotics of the  number of weighted partitions, when the Dirichlet generating function  for weights
 has a single  pole on the positive real axis.
Continuing \cite{GSE}, we derive asymptotics for the numbers of three basic types of decomposable combinatorial structures
(or, equivalently, ideal gas models in statistical mechanics) of size $n$, when their  Dirichlet generating functions have multiple simple poles on the positive real axis.
Examples to which our theorem applies include ones related to vector partitions and quantum field theory.
Our  asymptotic formula for the number of weighted partitions disproves  the belief accepted in the physics literature
that the main term in the asymptotics is determined by the rightmost pole.
\end{abstract}

\nin \section{Statement, history and motivation of the main result }
 \setcounter{equation}{0}
The goal of our work is the extension of Meinardus'  theorem to the case when the Dirichlet generating function for weights has more than one simple pole on the real axis. For this purpose we combine Meinardus' approach with the
probabilistic method of Khintchine, working in the framework of the  unified method of
asymptotic enumeration  of decomposable
combinatorial structures developed in \cite{GSE}.\\
The paper consists of two sections. In Section 1 we give a mathematical set up of the problem, the description of the probabilistic method used, the statement of our main result and a sketch of the relevant history. Section 2 is devoted to the proofs.

Our starting point is the following formalism.
Let $f$ be a
generating function of a nonnegative sequence $\{c_n,\ n\ge 0,\ c_0=1\}$:
\be
f(z)=\sum_{n\ge 0}c_n z^n.\la{sac}
\en
A specific feature of  structures considered is that
  $f$ has the following multiplicative form: \be
f=\prod_{k\ge 1}S_k,\la{sak1}\eu  where $S_k$ is a generating function for
some nonnegative sequence $\{d_k(j),\ j\ge 0,\ k\ge1\},$ i.e.\be
S_k (z)=\sum_{j\ge 0}d_k(j)z^{kj}, \ k\ge 1,\la{g3}\eu
where the sequence $\{d_k(j)\}$ is such that the infinite product $f$ in \refm[sak1] has a nonzero radius of convergence.
The above setting induces a sequence of probability measures $\mu_n, \ n\ge 1$ on the sequence of sets $\Omega_n, \ n\ge 1$ of integer partitions of $n$, such that $c_n$ is a partition function of $\mu_n$. Indeed,
denoting by $\eta_n=(j_1,\ldots,j_n)\in \Omega_n$ a partition of $n$ in the sense that $\sum_{k=1}^nkj_k=n$, we have from \refm[sac]-\refm[g3]
$$c_n=\sum_{\eta_n\in \Omega_n}\prod_{k=1}^nd_k(j_k).$$ Consequently, the corresponding measures $\mu_n$ are defined by
$$\mu_n(\eta_n)=c_n^{-1}\prod_{k=1}^nd_k(j_k).$$  The measures are called
multiplicative (Vershik \cite{V1}) or Gibbsian (Pitman \cite{Pi}).
 To explain the latter name we may view  $\mu_n$ as an equilibrium distribution of a physical system, presenting $\mu_n$ in the form
$$\mu(\eta_n)=c_n^{-1}\exp\big( \sum_{k=1}^n \log d_k(j_k)\big), \quad \eta_n\in\Omega_n,$$
so that $\sum_{j=1}^n \log d_k(j_k)$ is the total inner  energy of the system at the state $\eta_n.$
Since  the above expression for the total energy does not include  local potentials ascribed to interactions between groups of particles, the Gibbs measures $\mu_n$ characterize  systems without interactions.

 The measures $\mu_n$ are associated with a realm of models in statistical mechanics and in combinatorics.
The objective of the present paper is the asymptotics of $c_n$, as $n\to \infty.$ Our subsequent asymptotic analysis is based on the fundamental representation of $c_n$ for general multiplicative measures $\mu_n$:
\be
c_n=e^{n\delta} f_n(e^{-\delta})\PP\left(U_n=n\right),\quad n\ge
1,
 \la{rep}\eu
where $\delta>0$ is a free parameter,
\be
f_n=\prod_{k=1}^n S_k
\la{fndef}\eu
is the associated truncated generating function and, finally,
$$U_n:=\sum_{k=1}^{n}{Y_k},$$ where  $Y_k$ are independent integer-valued random variables distributed
in accordance with the given generating function $f$:
\be \PP(Y_k=jk)= \frac{d_k(j) e^{-\delta kj}}{S_k(e^{-\delta })},
\quad j\ge 0,\quad k\ge 1.\la{Y}\eu

The representation \refm[rep] which was  derived in \cite{GSE} is an analogue of Khintchine's representation of partition functions $c_n$  for particular models of statistical mechanics. In accordance with  Khintchine's probabilistic method,
 we will choose the free parameter $\delta=\delta_n$ to be the solution
of the equation
$$
\BE U_n=n,\ n\ge 1,
$$
which, by virtue of the above setting,  can be written as \be -\Big(\log(f_n(e^{-\delta}))\Big)^\prime_{\delta=\delta_n}=n,\quad n\ge 1,\la{der} \eu
where the derivative is taken with respect to $\delta$.
By the laws of thermodynamics, the function $\Phi(\delta):=\log\Big(e^{n\delta}f_n(e^{-\delta})\Big)$
is the entropy of the physical system considered. It is simple to show that  $\Phi^{\prime\prime}(\delta)>0,\ \delta>0$ which tells us that the solution $\delta_n$ of \refm[der] is the point of minimum of the entropy $\Phi(\delta)$.
This important meaning of the aforementioned choice of the free parameter was revealed by Khintchine (\cite{Kh}, Ch.6), for simple models of statistical mechanics,
allowing him to provide a probabilistic derivation of the second law of thermodynamics.

Our study is restricted to  three classic generating functions $f^{(i)},\ i=1,2,3$
associated with three types of  decomposable combinatorial structures, each of which can be interpreted as a model of ideal gas from statistical mechanics. These structures are weighted partitions (Bose Einstein statistics), selections (Fermi-Dirac statistics) and assemblies (Maxwell-Boltzmann statistics), respectively. Defining the functions $\nicef^{(i)}(\delta)=f^{(i)}\left(e^{-\delta}\right), \ \delta>0$, for $i=1,2,3$, we have
\begin{eqnarray}
\nicef^{(1)}(\delta)&=&\prod_{k\ge 1}{(1-e^{-k\delta})^{-b_k}},\la{F1}\\
\nicef^{(2)}(\delta)&=&\prod_{k\ge 1}{(1+e^{-k\delta})^{b_k}},\la{F2}\\
\nicef^{(3)}(\delta)&=&\exp\left(\sum_{k\ge
1}b_ke^{-k\delta}\right).\la{F3}
\end{eqnarray}
 Clearly, by  definition \refm[fndef], the quantities   $f_n^{(i)}(e^{-\delta})$, $i=1,2,3$ are given by the expressions  \refm[F1]-\refm[F3],  with products  restricted
to the range $1\leq k\leq n$.
Correspondingly,   the distributions for the random variables $\frac{1}{k}Y_k$ in \refm[rep] are of the following
three types: Negative
 Binomial $\left(b_k;e^{-\delta k}\right),$
Binomial $\left(b_k;\frac{e^{-\delta k}}{1+e^{-\delta k}}\right)$
and Poisson $\left(b_ke^{-\delta k}\right),$
and the equation \refm[der] takes the forms:
\begin{eqnarray}
\sum_{k=1}^{n}{\frac{kb_ke^{-k\delta^{(1)}_n}}{1-e^{-k\delta^{(1)}_n}}}
&=&n, \label{def1deltan}\\
 \sum_{k=1}^{n}{\frac{kb_ke^{-\delta^{(2)}_n
k}}{1+e^{-\delta^{(2)}_n k}}}&=&n,  \label{def2deltan}\\
\sum_{k=1}^{n} {kb_ke^{-\delta^{(3)}_n k}}&=&n.
\label{def3deltan}
\end{eqnarray}
A common feature of the three models considered is that each of them is given by one sequence of parameters, which is $\{b_k\ge 0,\ k\ge 1\}$.
In the context of statistical mechanics $b_k$ is the weight of the  energy level $k,$
while in combinatorics $b_k$ is the weight prescribed to a  indecomposable component of size $k$.

 Following Meinardus' approach we define two generating functions  for the sequence $\{b_k\}$: the Dirichlet  generating function $D$ and the power generating function G, given by \be D(s)=\sum_{k=1}^\infty b_k
k^{-s}, \quad s=\sigma+it \label{D},\eu
\be G(z)= \sum_{k=1}^\infty b_kz^k,\quad \vert z\vert<1.   \la{fak}\eu
We assume that $D$ and $G$ satisfy the conditions $(I)- (III)$ below:

{\bf Condition $(I)$.}
The Dirichlet generating function $D(s), \ s=\sigma+it$ is analytic in the half-plane $\sigma>\rho_r>0$ and
it  has $r\geq 1$ simple poles at positions $0<\rho_1<\rho_2<\ldots<\rho_r,$ with positive residues
$A_1,A_2,\ldots,A_r$ respectively. Moreover,
 there is a constant $0<C_0\le 1,$ such
that the function $D(s)$, $s=\sigma+it$, has a meromorphic continuation
 to the half-plane \be\label{Hdef}
\calh=\{s:\sigma\geq-C_0\} \eu on which it is analytic except for
the above $r$ simple poles.

{\bf Condition $(II)$.}
 There is a constant $C_1>0$ such that
\be\label{imagbound} D(s)=O\left(|t|^{C_1}\right),\quad t\to\infty
\eu uniformly for $s=\sigma+it\in\calh$.

{\bf Condition $(III).$}
For $\delta>0$ small enough and any
$\epsilon>0,$
$$
2\sum_{k=1}^\infty
 {b_k{e^{-k\delta}\sin^2(\pi k\alpha)}}\ge\left(1+{\rho_r\over 2}+\epsilon\right)
{\cal M}^{(i)} |\log \delta|,$$
$$ \quad  \sqrt{\delta} \leq
|\alpha|\leq 1/2,\quad i=1,2,3,
$$
where  the constants ${\cal M}^{(i)}$ are defined by
$${\cal M}^{(i)}=
  \begin{cases}
    {4\over{\log 5}} , & \text{if } \quad i=1,\\
    4 , & \text{if}\quad i=2,\\
     1, & \text{if} \quad i=3.
  \end{cases}
$$

{\bf Remark.} The main difference between conditions $(I)-(III)$ and the conditions used in the original paper of Meinardus \cite{M} is that condition $(I)$ allows $D(s)$ to have multiple poles on the positive real axis. Another difference between our conditions and the conditions used by Meinardus is condition $(III)$, which as was shown in \cite{GSE},
 is  weaker than the corresponding condition in \cite{M}.

To formulate our main result we need some more notations. Define the finite set

\be\la{sd1} \tilde{\Upsilon}_r=
\left\{\sum_{k=0}^{r-1} d_k(\rho_r-\rho_k):\   d_k\in\BZ_+,\ \sum_{k=0}^{r-1}d_k\geq 2\right\}\cap \big(0,\rho_r+1\big],
\en
where we have set $\rho_0=0$ and let  $\BZ_+$ denote the set of nonnegative integers. Let $0<\alpha_1<\alpha_2<\ldots<\alpha_{|\tilde{\Upsilon}_r|}\le \rho_r+1$ be all ordered numbers forming the
set $\tilde{\Upsilon}_r$. Clearly, $\alpha_1=2(\rho_r-\rho_{r-1}),$ if the set $\tilde{\Upsilon}_r$ is not empty.
We also
 define  the finite set
\be\la{sd2} \Upsilon_r=\tilde{\Upsilon}_r\cup\{\rho_r-\rho_k:\ k=0,1,\ldots,r-1\},
\en
observing that  some of the differences $\rho_r-\rho_k,\ k=0,\dots,r-1$ may fall into the set $\tilde{\Upsilon}_r$. We let $0<\lambda_1<\lambda_2<\ldots<\lambda_{|\Upsilon_r|}$
be all ordered numbers forming the
set $\Upsilon_r$. In addition to its appearance in (\ref{sd2}),
the set  $\tilde{\Upsilon}_r$ plays an auxiliary role in the proof of
Proposition~\ref{deltaexp} below.

Now that we have defined the quantities $\lambda_s$, we are ready to state our main result.
\begin{theorem}({\bf Meinardus type theorem with multiple singularities}).
\label{main}

Suppose that the conditions (I)-(III) hold. Then the following asymptotic formulae for $c_n^{(i)},\ i=1,2,3$
are valid:
$$
c_n^{(i)}\sim H^{(i)}n^{-\frac{2+\rho_r-2D(0){\bf 1}(i)}{2(\rho_r+1)}}\exp \Big(\sum_{l=0}^{r}P_l^{(i)}n^{\frac{\rho_l}{\rho_r+1}}+
\sum_{l=0}^r h_l^{(i)}\sum_{s:\lambda_s\le \rho_l} K^{(i)}_{s,l} n^{\frac{\rho_l-\lambda_s}{\rho_r+1}}\Big),
$$
\be n\to \infty, \la{xhu}\en
where $${\bf 1}(i)=\left\{
                                                                              \begin{array}{ll}
                                                                                1, & \hbox{i=1}; \\
                                                                                0, & \hbox{i=2,3},
                                                                              \end{array}
                                                                            \right.
$$
where the powers $\frac{\rho_l-\lambda_s}{\rho_r+1}$ do not depend on $n$
and where the coefficients $H^{(i)}, P_l^{(i)}, K^{(i)}_{s,l},h_l^{(i)}$ that do not depend on $n$ either,  are implicitly defined in the course of the proof of the theorem in Section 2 below.

\end{theorem}
{\bf A brief history and motivation.}
 The famous Meinardus' theorem  published in 1954, in \cite{M}, provided an asymptotic formula for $c_n^{(1)}$, under three assumptions, call them $(I^\prime)-(III^\prime)$, which are prototypes of our conditions $(I)-(III).$
The assumption $(I^\prime)$ differs from $(I)$ in that it requires that $D$ has only one simple positive pole, the assumption
$(II^\prime)$ is identical to $(II)$, while  the technical condition  $(III)$ is a weaker form of $(III^\prime)$.
The asymptotic formula  of \cite{M} preceded by the great Hardy-Ramanujan exact asymptotic expansion (1918) for the number of integer partitions, was widely discussed in the literature on enumerative combinatorics and statistical mechanics.

 A formal motivation for our extension of the original Meinardus' theorem to the case when $D$ has many poles
comes from the  basic restriction of the theorem, implied by its  condition $(I^\prime)$. Due to this condition the theorem  being valid  for  the weight sequences of the form $b_k=ak^{r-1},\ a,r>0,$  fails for linear combinations of such $b_k$'s, because the induced
functions $D$ have many  poles.

In modern physics, the need of such extension arises in the  counting of BPS operators, a line of research that originated in field theory in the 2000's (see \cite{LR}).
In the language of combinatorics, problems treated in this context (see \cite{BFHH,LR}) belong to weighted partitions with weights of the form $b_k={k+l\choose l} $, for some $l\ge 2.$
 The Dirichlet generating function of such a weighted partition has $r=l+1$ simple poles at integer points $1,2,\ldots,l+1$.

In the physics literature  it has been accepted the belief
that in the case of $r>1$ poles $0<\rho_1,\ldots,\rho_r$ the number  of weighted partitions $c_n^{(1)}=c_n^{(1)}(\rho_1,\ldots,\rho_r)$ is asymptotically equivalent to the number $c_n^{(1)}(\rho_r)$ of partitions corresponding to the model with one rightmost pole $\rho_r$ (see (6.9) in \cite{BFHH} and (5.19) in \cite{LR}). Our asymptotic formula \refm[xhu] for $c_n^{(i)}, \ i=1,2,3$   disproves this belief for all three types of the structures considered.
Nevertheless, the basic idea behind this belief is correct, in the sense that
\be \log c_n^{(i)}(\rho_1,\ldots,\rho_r)
\sim \log c_n^{(i)}(\rho_r)
\sim
P_r^{(i)}
n^{\frac{\rho_r}{\rho_r+1}}, \ n\to \infty,\la{gad}
\en
by virtue of \refm[xhu].
The proof of Theorem~\ref{main} shows that
$$
P_r^{(i)}>0, \ i=1,2,3.
$$

Finally, we note that an interesting extension of Meinardus' theorem in a related but a quite different direction was obtained
in \cite{MW}\ and \cite{TT}. In both papers it is studied  the problem of asymptotic enumeration of the number of unweighted partitions ($b_k=1,\ k\ge 1)$ of $n$
with summands belonging to a given infinite set of positive integers $\Lambda=\{\lambda_1,\ldots,\lambda_k,\ldots\}$ and with associated spectral zeta function $D(s)=\sum_{l\ge 1} \lambda_l^{-s}$. In \cite{MW} the set $\Lambda$ represents numbers with missing digits, so that the generating function  $D$ has
 simple  equidistant poles $\alpha+2\pi ik\omega, \ k\in {\cal Z}$,
for a fixed $\omega>0$, on the line $\Re( s)=\alpha$ for some $\alpha>0,$ and it is analytic in the half-plane $\Re( s)>\alpha$. In \cite{TT} motivated by the enumeration problems in the setting of quantum mechanics, the function $D$ has simple real poles at the integers $n,n-1,\ldots,1,-1,\ldots$.
\nin \section{Proofs.}
The most difficult part of the proof is contained in the Subsection 2.2 in which we establish the asymptotic formulae for the solutions $\delta_n^{(i)},\ i=1,2,3$ in the case of $r\ge 1$ simple poles.

 In what follows we always suppose that conditions $(I)-(III)$ are satisfied.

\nin \subsection{Asymptotics of generating functions, as $\delta\to 0^+.$}

\begin{lemma} \label{estimates}
\noindent  (i) As $\delta\to 0^+$,
\be \label{prod}
\nicef^{(i)}(\delta)=\exp\left(\sum_{l=0}^r h_l^{(i)}\delta^{-\rho_l}-\big(D(0)\log\delta\big){\bf 1}(i)+M^{(i)}(\delta;C_0)\right),\en
where
$$
\rho_0=0,\quad {\bf 1}(i)=\left\{
                                                                              \begin{array}{ll}
                                                                                1, & \hbox{i=1}; \\
                                                                                0, & \hbox{i=2,3},
                                                                              \end{array}
                                                                            \right.
$$
$$h_l^{(i)}=\left\{
              \begin{array}{ll}

A_l\Gamma(\rho_l)\zeta(\rho_l+1), & \hbox{i=1, l=1,\ldots,r;} \\
A_l\Gamma(\rho_l)(1-2^{-\rho_l})\zeta(\rho_l+1), & \hbox{i=2, l=1,\ldots,r;} \\
                A_l\Gamma(\rho_l), & \hbox{i=3, l=1,\ldots,r},
              \end{array}
            \right.
$$
$$h_0^{(i)}=\left\{
              \begin{array}{ll}
                D^\prime(0), & \hbox{i=1}; \\
                D(0)\log2, & \hbox{i=2}; \\
                D(0), & \hbox{i=3}
              \end{array}
            \right.
$$
and
\begin{eqnarray} M^{(1)}(\delta;C_0)&=&\frac{1}{2\pi
i}\int_{-C_0-i\infty}^{-C_0 +i\infty}
\delta^{-s}\Gamma(s)\zeta(s+1)D(s)ds, \la{Mexp1}\\
M^{(2)}(\delta;C_0)&=&\frac{1}{2\pi
i}\int_{-C_0-i\infty}^{-C_0+i\infty}
\delta^{-s}\Gamma(s)(1-2^{-s})\zeta(s+1)D(s)ds,\la{Mexp2}\\
M^{(3)}(\delta;C_0)&=&\frac{1}{2\pi
i}\int_{-C_0-i\infty}^{-C_0+i\infty} \delta^{-s}\Gamma(s)D(s)ds, \la{Mexp3}
\end{eqnarray}
with $M^{(i)}(\delta;C_0) =O(\delta^{C_0}),\ i=1,2,3.$

(ii) The  asymptotic expressions for the derivatives
$$\Big(\log\nicef^{(i)}(\delta)\Big)^{(k)}, \ i=1,2,3,\ k=1,2,3$$
are given by the formal differentiation of the logarithm of
\refm[prod], with \\  $(M^{(i)}(\delta;C_0))^{(k)}=O(\delta^{C_0-k}),\ i=1,2,3,\ k=1,2,3.$
\end{lemma}
\proof
 The proof of the claim $(i)$
which  we only sketch is  similar to the proof of the first part  of Lemma 2 in
\cite{GSE}.
 Following the Meinardus approach, we will use the fact that
$e^{-u}$, $u>0$, is the Mellin transform of the Gamma function:
\be e^{-u}=\frac{1}{2\pi i}\int_{v-i\infty}^{v+i\infty}
u^{-s}\Gamma(s)\,ds,\quad u>0,\ \Re(s)=v>0. \la{Mellin} \en
Expanding $\log(1-e^{-\delta k})$ and $\log(1+e^{-\delta k}) $ from \refm[F1] and \refm[F2], respectively,
in terms of $e^{-\delta k}$ and then substituting \refm[Mellin] in \refm[F1] - \refm[F3],
with $v=\rho_r+\epsilon,$ with any $\epsilon>0,$
gives the desired integral representations of the functions
$\log
~\nicef^{(i)}(\delta), \quad i=1,2,3$ for all $\delta>0$:
\begin{eqnarray} \log ~\nicef^{(1)}(\delta)&=&\frac{1}{2\pi
i}\int_{\epsilon +\rho_r-i\infty}^{\epsilon +\rho_r+i\infty}
\delta^{-s}\Gamma(s)\zeta(s+1)D(s)ds, \la{intrep}\\
\log ~\nicef^{(2)}(\delta)&=&\frac{1}{2\pi
i}\int_{\epsilon+\rho_r-i\infty}^{\epsilon + \rho_r+i\infty}
\delta^{-s}\Gamma(s)(1-2^{-s})\zeta(s+1)D(s)ds,\la{intrep2}\\
\log ~\nicef^{(3)}(\delta)&=&\frac{1}{2\pi
i}\int_{\epsilon+\rho_r-i\infty}^{\epsilon+\rho_r+i\infty} \delta^{-s}\Gamma(s)D(s)ds. \la{intrep3}
\end{eqnarray}
Next, we apply the residue theorem for the above integrals, in the complex domain $ -C_0\le \Re( s)\le \rho_r+\epsilon, $ with $0\le
C_0<1, \epsilon>0.$
By virtue of  condition $(I)$,
 the integrands in \refm[intrep]-\refm[intrep3] have in the aforementioned domain $r$
simple poles at $\rho_l>0, \ l=1,\ldots,r$.    The corresponding residues at $s=\rho_l$ in the three cases are equal to:

\begin{itemize}
\item[] $A_l\delta^{-\rho_l}\Gamma(\rho_l)\zeta(\rho_l+1),$ $l=1,\ldots,r,$ {\rm \ for \ } $i=1$,

\item[] $A_l\delta^{-\rho_l}\Gamma(\rho_l)(1-2^{-\rho_l})\zeta(\rho_l+1),$ $l=1,\ldots,r,$ {\rm \ for \ } $i=2$,

\item[] $A_l\delta^{-\rho_l}\Gamma(\rho_l),$ $l=1,\ldots,r,$ {\rm \ for \ } $i=3$.
\end{itemize}

Besides,
from the Laurent expansions at $s=0$ of the Riemann Zeta
 function $\zeta(s+1)=\frac{1}{s}+\gamma+\ldots$ and
the Gamma function $\Gamma(s)=\frac{1}{s}-\gamma+\ldots,$ where
$\gamma$ is
 Euler's constant,  the integrands have  also a pole at $s=0$, which is a simple one in the cases $i=2,3$ and is of  a second order
  in the case $i=1$. The residues at $s=0$ are equal to

\begin{itemize}
\item[] $D^\prime(0)-D(0)\log\delta$ {\rm \ for \ } $i=1$,

\item[] $D(0)\log2$ {\rm \ for \ } $i=2$,

\item[] $D(0)$ {\rm \ for \ } $i=3$.
\end{itemize}

Finally, applying condition $(II)$ shows that in all three cases the integrals of the above integrands, over the horizontal contour $-C_0\le\Re(s)\le \epsilon+\rho_r$, $\Im(s)=t$, tend to zero, as $t\to \infty,$ for
any fixed $\delta>0.$   This gives the claimed formulae \refm[prod] for $i=1,2,3$, where the remainder terms $M^{(i)}$ are
integrals of the aforementioned integrands,  taken over the vertical contour $-C_0+it,\ -\infty<t<\infty$.

The proof of $(ii)$ proceeds as in the proof of Lemma 2 in \cite{GSE}.
\hfill\qed

The asymptotics of the generating functions for weighted partitions has almost a century long history  of its own which is enlightened in the  remarkable monograph \cite{K}, p.228-229,  by Korevaar, (see also \cite{GSE}).
\nin \subsection{Asymptotics of $\delta_n$. }
{\bf Preliminaries.} First we note that for $i=1,2,3$ equations
\refm[def1deltan] - \refm[def3deltan] have unique solutions $\delta_n^{(i)}$,
such that $\delta_n^{(i)}\to 0$ as $n\to \infty$. The proof
of this is the same as  the proof of part $(ii)$ of Lemma~2 in \cite{GSE}.
We call any $\tilde{\delta}_n,$ such that \be \big(-\log f_n(e^{-\delta})\big)_{\delta=\tilde{\delta}_n}^\prime -n\to 0,\quad n\to \infty \la{asympt}\en
an asymptotic solution of \refm[der].
We demonstrate that in all three cases considered  it is sufficient for \refm[asympt] that $\tilde{\delta}_n$ obeys the condition
\be (- \log \nicef(\delta))_{\delta=\tilde{\delta}_n}^\prime -n\to 0, \ n\to \infty. \la{aut} \en
In fact, by Lemma~\ref{estimates},$$ \big(- \log \nicef^{(i)}(\delta)\big)^\prime\sim h_r^{(i)}\rho_r \delta^{-\rho_r-1},\quad \delta\to 0,\quad i=1,2,3,$$ from which it follows that \refm[aut] implies
\be\tilde{\delta}_n^{(i)}\sim (h_r^{(i)}\rho_r)^{\frac{1}{\rho_r+1}}n^{-\frac{1}{\rho_r+1}}, \ n\to \infty, \quad i=1,2,3.\la{caf}\en
We have,  for an arbitrary multiplicative measure and all $\delta>0,$ \be\log f_n(e^{-\delta})=\log \nicef(\delta)-\sum_{k=n+1}^\infty \log S_k(\delta).\la{ser}\en
By using the Wiener-Ikehara theorem as in the proof of Theorem~2 of \cite{GSE},
one may show that in our setting  $b_k=o(k^{\rho_r})$. (In this connection note that for the application of the theorem only the rightmost pole matters). Using this fact, we obtain
 from \refm[F1] for the case $i=1$:
\begin{eqnarray}
\sum_{k=n+1}^\infty\Big(- \log S_k^{(1)}(\tilde{\delta}^{(1)}_n)\Big)^\prime&=&
\sum_{k=n+1}^\infty{\frac{kb_ke^{-k\tilde{\delta}^{(1)}_n}}{1-e^{-k\tilde{\delta}^{(1)}_n}}}\nonumber\\&=&
 o\left(\int_{n+1}^\infty {\frac{u^{\rho_r+1}e^{-u\tilde{\delta}^{(1)}_n}}{1-e^{-u\tilde{\delta}^{(1)}_n}}}du\right)\nonumber\\&=&
\Big(\tilde{\delta}^{(1)}_n\Big)^{-\rho_r-2}o\left(\int_{(n+1)\tilde{\delta}^{(1)}_n}^\infty
\frac{w^{\rho_r+1}e^{-w}}{1-e^{-w}}\,dw\right)\nonumber\\
&=&o(1),\quad n\to \infty,
\la{arg}
\end{eqnarray}
by \refm[caf],
since $n\tilde{\delta}^{(1)}_n\to \infty$. This together with \refm[ser] say that \refm[aut] implies
\refm[asympt] in the case $i=1$. For the cases $i=2,3$ the proof of the claim is similar.\\
As in \cite{GSE}, we will use \refm[aut] to derive asymptotic expansions for $\tilde{\delta}^{(i)}_n,\ i=1,2,3$
in our setting.
It is necessary to establish that the error of   approximating the exact solution $\delta_n$ by the asymptotic solution $\tilde{\delta}_n$ is of order $o(n^{-1}).$
By the monotonicity in $\delta>0$ (for any fixed $n\ge 1$) of the functions $\log f_n^{(i)}(e^{-\delta}), \ i=1,2,3$ and their derivatives,
 it follows from \refm[caf] that
\be \delta_n^{(i)}\sim \tilde{\delta}_n^{(i)},\ n\to \infty,\ i=1,2,3. \la{gorn}\en
Applying the Mean Value Theorem, we obtain
\be \left| \Big(\log f_n\big(e^{-\delta_n^{(i)}}\big)\Big)^\prime-\Big(\log f_n\big(e^{-\tilde{\delta}_n^{(i)}}\big)\Big)^\prime\right|=
\left|(\delta_n^{(i)}-\tilde{\delta}_n^{(i)}) \Big(\log f_n(e^{-u_n})\Big)^{\prime\prime}\right|,\la{ghj}\en
where
$$u_n\in[\min(\delta_n^{(i)},\tilde{\delta}_n^{(i)}),\max(\delta_n^{(i)},\tilde{\delta}_n^{(i)})].$$ By virtue of the definitions \refm[der], \refm[asympt], the left hand side of \refm[ghj] tends to $0$, as $n\to \infty$.
It follows from the preceding arguments and Lemma~\ref{estimates} that for all three cases,
\be\Big(\log f_n(e^{-u_n})\Big)^{\prime\prime}\sim \rho_r(\rho_r+1)h_r^{(i)}(\delta_n^{(i)})^{-\rho_r-2}= O(n^{\frac{\rho_r+2}{\rho_r+1}}).\la{second}\eu
Combining \refm[ghj] with \refm[second], gives the desired estimate
\be \left|\delta_n^{(i)}-\tilde{\delta}^{(i)}_n\right|=o(n^{-1}),\quad i=1,2,3. \la{nmi}\en
 {\bf The asymptotic solution.} Now our efforts will be devoted to derive  the asymptotic expansions for $\tilde{\delta}_n^{(i)}, \ i=1,2,3$ which by \refm[nmi], will provide asymptotic expansions of $\delta_n^{(i)}$ up to order $o(n^{-1}).$
For the sake of brevity, we will use in the course of the proof the notations $$\hat{h}_l^{(i)}=\rho_lh_l^{(i)}, \quad l=1,\ldots r,\quad i=1,2,3,$$
and
$$\ \hat{h}_0^{(i)}=\left\{
                                           \begin{array}{ll}
                                             D(0), & \hbox{i=1};\\
                                             0, & \hbox{i=2,3.}
                                           \end{array}
                                         \right.
$$

\begin{proposition}\label{deltaexp}  Suppose that the sequence $b_k\geq 0$, $k\geq 1$ is such that the
associated Dirichlet generating function $D$ satisfies the conditions
$(I)$ and $(II)$ of Theorem~\ref{main}. Then,
the asymptotic expansions of the $\delta_n^{(i)}$ up to terms of order $o(n^{-1})$ are given by
\be
\delta_n^{(i)}=\big(\hat{ h}_r^{(i)}\big)^{\frac{1}{\rho_r+1}}n^{-\frac{1}{\rho_r+1}}+ \sum_{s=1}^{\vert\Upsilon_r\vert} K^{(i)}_s n^{-\frac{1+\lambda_s}{\rho_r+1}}+o(n^{-1}), \ i=1,2,3, \la{sd6}
\en
where  $K^{(i)}_s$ and $\lambda_s$ do not depend on $n$, and the powers $\lambda_s$ are as defined in \refm[sd2].

\end{proposition}
\proof   By $(ii)$ of Lemma~\ref{estimates} we have
$$\Big(-\log\nicef^{(i)}(\delta)\Big)^{\prime}= \sum_{l=0}^r \hat{h}_l^{(i)}\delta^{-\rho_l-1}+ \Big(M^{(i)}(\delta;C_0)\Big)^\prime,\ i=1,2,3.$$
(Here and in what follows $\Big(M^{(i)}(\delta;C_0)\Big)^\prime$ denotes the derivative with respect to $\delta$).
By \refm[caf], the condition \refm[aut] may be  rewritten as
\begin{eqnarray} \Delta_n^{(i)}&:=&
\frac{n(\tilde{\delta}_n^{(i)})^{\rho_r+1}-\sum_{k=0}^r
\hat{h}_k^{(i)}(\tilde{\delta}_n^{(i)})^{\rho_r-\rho_k}-
(\tilde{\delta}_n^{(i)})^{\rho_r+1}\Big(M^{(i)}(\tilde{\delta}_n^{(i)};C_0)\Big)^\prime}{(\tilde{\delta}^{(i)}_n)^{\rho_r+1}}\nonumber\\&\sim&
\frac{n(\tilde{\delta}_n^{(i)})^{\rho_r+1}-\sum_{k=0}^r
\hat{h}_k(\tilde{\delta}_n^{(i)})^{\rho_r-\rho_k}-(\tilde{\delta}_n^{(i)})^{\rho_r+1}
\Big(M^{(i)}(\tilde{\delta}_n^{(i)};C_0)\Big)^\prime}{\hat{h}_r
n^{-1}}\nonumber\\
&=&o(1),\quad n\to \infty, \quad i=1,2,3,  \label{Dell}
\end{eqnarray}
which is equivalent to saying that the numerator in \refm[Dell] is $o(n^{-1}).$ In the course of the proof we will set
\be
\tilde{\delta}_n^{(i)}=Q_n^{(i)}+u_n^{(i)}, \ i=1,2,3,
\la{Qplusu}\eu for some
$Q_n^{(i)}\sim
\big(\hat{ h}_r^{(i)}\big)^{\frac{1}{\rho_r+1}}n^{-\frac{1}{\rho_r+1}}\sim \tilde{\delta}_n^{(i)}$ and $u_n^{(i)}=o(\tilde{\delta}_n^{(i)})=o(Q_n^{(i)}), \ i=1,2,3$
to be determined.
 We will then set
$Z_n^{(i)}=\frac{u_n^{(i)}}{Q_n^{(i)}}$. Since the rest of the proof goes the same way for all three cases considered, we suppress  the index $(i)$ in the forthcoming expansions.
To motivate the choice of $Q_n^{(i)}$ and $u_n^{(i)}$, we observe that
the insertion of $\tilde{\delta}_n^{(i)}$ of the form
\refm[Qplusu] in \refm[Dell], and the use  the binomial expansion gives
\begin{eqnarray}
\hat{h}_r n^{-1}\Delta_n&\sim &
nQ_n^{\rho_r+1}-\sum_{k=0}^r \hat{h}_kQ_n^{\rho_r-\rho_k} \non\\
&& +\sum_{m\ge 1}\binom{\rho_r+1}{m} (nQ_n^{\rho_r+1})Z_n^m- \la{wria}\\
&& -\sum_{k=0}^{r-1}
\hat{h}_kQ_n^{\rho_r-\rho_k}\Big(\sum_{m\ge 1}\binom{\rho_r-\rho_k}{m}Z_n^m\Big)-\la{wrib} \\
&&-\Big(M(\tilde{\delta}_n;C_0)\Big)^\prime Q_n^{\rho_r+1}\sum_{m\ge
0}\binom{\rho_r+1}{m}Z_n^m\la{wric}\\&=& nQ_n^{\rho_r+1}-\sum_{k=0}^r \hat{h}_kQ_n^{\rho_r-\rho_k}+ L_n,\la{wri}
\end{eqnarray}
where $L_n$ denotes \refm[wria] + \refm[wrib] + \refm[wric].
Observe that if $Q_n$ has been chosen, then $L_n$ depends on $Z_n$ alone.

Our plan for the remainder of the proof consists of the following three steps:

Step 1. Determine a $Q_n\sim\delta_n$ satisfying  the condition
\be nQ_n^{\rho_r+1}-\sum_{k=0}^r
\hat{h}_kQ_n^{\rho_r-\rho_k}=o(n^{-1}).\la{condd}\en

Step 2.  Show that, under $Q_n$ obeying \refm[condd], one may choose  $Z_n,$  so that $L_n$ defined in \refm[wri] satisfies  $L_n =o(n^{-1}).$

Step 1 and Step 2 construct an asymptotic solution  $\tilde{\delta}^{(i)}_n$ of the form \refm[Qplusu].
By \refm[nmi], this  provides the asymptotic expansion for the solutions $\delta^{(i)}_n, \ i=1,2,3$ up to the of order $o(n^{-1}).$
At this point it is convenient to introduce a variable
$$z=z_n=n^{-\frac{1}{\rho_r+1}}.$$

Step 3. Express the obtained asymptotic formulae for $\delta^{(i)}_n, \ i=1,2,3$ as  finite linear combinations of powers $z$.
\vskip .5cm
{\bf Step 1}.
For the sake of convenience, we introduce another notation. Given a function $Q_n(z)$ of $z$,
we define
\be\tilde{Q}(z)=z^{-1}Q_n(z).\la{tilQdef}\eu
Then \refm[condd] can be written
as \be \tilde{Q}^{\rho_r+1}(z)-
\sum_{k=0}^r\hat{h}_k\tilde{Q}^{\rho_r-\rho_{k}}
(z)z^{\rho_r-\rho_{k}}=o(n^{-1}). \la{cond}\en It follows
from \refm[cond] that
 $\tilde{Q}^{\rho_r+1}(0)=\hat{h}_r,$  while \refm[tilQdef] implies
$\tilde{Q}^{\rho_r+1}(z)=O(1)$. We will
verify that \refm[cond] is satisfied with $\tilde{Q}$   given by
 \be
\tilde{Q}^{\rho_r+1}(z)=\tilde{Q}^{\rho_r+1}(0)
+\sum_{k=0}^{r-1}\hat{h}_k
\tilde{Q}^{\rho_r-\rho_k}(0)z^{\rho_r-\rho_k}+V(z),\la{hal}\en
where $V(z)\equiv 0,$ if $\tilde{\Upsilon}_r=\varnothing,$
but where otherwise $V(z)$ is of the form
 \be V(z)=  \sum_{m=1}^{\vert\tilde{\Upsilon}_r\vert} B_m z^{\alpha_m},\la{jur}\en where
the coefficients $B_m$ do not depend on $z$, and where   the powers
$\alpha_1<\alpha_2<\cdots<\alpha_{|\tilde{\Upsilon}_r|}$ comprise the set $\tilde{\Upsilon}_r$ defined by \refm[sd1].
 Denote
$$P(z)=\sum_{k=0}^{r-1}\hat{h}_k
\tilde{Q}^{\rho_r-\rho_k}(0)z^{\rho_r-\rho_k}$$ and assume that
$\tilde{Q}(z)$ is defined by
\refm[hal] for some function $V(z)$ of the form \refm[jur] with undetermined coefficients $B_m$.  Then the  condition  \refm[cond] can be
 written as
\be V(z)-\sum_{k=0}^{r-1}\hat{h}_k
z^{\rho_r-\rho_k}\Big(\tilde{Q}^{\rho_r-\rho_k}(z)-\tilde{Q}^{\rho_r-\rho_k}(0)\Big)=o(n^{-1}).\la{con1}\en
Moreover, we have
  \begin{eqnarray}
\tilde{Q}^{\rho_r-\rho_k}(z)-\tilde{Q}^{\rho_r-\rho_k}(0)&=&\tilde{Q}^{\rho_r-\rho_k}(0)\Bigg(
\Big(1+\frac{P(z)+V(z)}{\tilde{Q}^{\rho_r+1}(0)}\Big)^{\frac{\rho_r-\rho_k}{\rho_r+1}}-1\Bigg)\nonumber\\&=&
  O(z^{\rho_r-\rho_{r-1}}),\quad k=0,\ldots,r-1,                      \la{con2}\end{eqnarray}
where the last asymptotics follows from the definitions of $P(z),V(z)$. Hence \refm[con1] implies
  $V(z)=O(z^{2(\rho_r-\rho_{r-1}})),$ in accordance with the fact that $\alpha_1=2(\rho_r-\rho_{r-1}),$
provided the set $\tilde{\Upsilon}_r$ is not empty.


It follows from the definition of the set $\tilde{\Upsilon}_r$ that if we can find $V(z)$ of the form \refm[jur] such that
$$B_m-[z^{\alpha_m}]\Bigg(\sum_{k=0}^{r-1}\hat{h}_k
z^{\rho_r-\rho_k}\Big(\tilde{Q}^{\rho_r-\rho_k}(z)-\tilde{Q}^{\rho_r-\rho_k}(0)\Big)\Bigg)=0,$$\be
m=1,\ldots,\vert \tilde{\Upsilon}_r\vert,\la{ec}\en
where  the second term in \refm[ec] is the coefficient of $z^{\alpha_m}$
in the expression in parentheses,
then we will have found the required  $V(z)$ satisfying  \refm[con1].
Using \refm[con2], the equation \refm[ec] with $m=1$, $\alpha_1=2(\rho_r-\rho_{r-1})$ can be written as
$$B_1=\hat{h}^2_{r-1}\frac{(\rho_r-\rho_{r-1})\tilde{Q}^{\rho_r-2\rho_{r-1}-1}(0)}{\rho_r+1}.$$
The crucial fact is that, by virtue of  \refm[con2],  the coefficient
of $z^{\alpha_m}$ in  \refm[ec] does
not depend on $B_s$ for $s\geq m$, so that
$$B_m+ q_m(B_1,\ldots, B_{m-1})=0,\ \
m=2,\ldots,\vert\tilde{\Upsilon}_r\vert,$$ where $q_m$ are polynomials of $B_1,\ldots, B_{m-1}$.
Thus, the $B_m$ can be defined recursively from \refm[ec], resulting in a function $V(z)$ of the form \refm[jur] satisfying \refm[con1].

Finally, with this choice of $V(z)$, we obtain
from \refm[con2] $$ \tilde{Q}(z)=\tilde{Q}(0)
\left(1+\frac{P(z)+V(z)}{\tilde{Q}^{\rho_r+1}(0)}\right)^{\frac{1}{
\rho_r+1}}.
 $$
Recalling \refm[tilQdef], this gives \be
Q_n=\hat{h}_r^{\frac{1}{\rho_r+1}}z\Big(1+\frac{P(z)+V(z)}{\hat{h}_r}\Big)^{\frac{1}{
\rho_r+1}}, \la{soft} \en
which, by our construction, satisfies \refm[condd].

{\bf Step 2.} We first need several definitions.  Under $Q_n$ given by \refm[soft], define
\be C_n(m)= \binom{\rho_r+1}{m}nQ_n^{\rho_r+1}\sim {{\rho_r+1}\choose m}
\hat{h}_r\neq 0, \quad m\ge 1, \la{Cdef}\eu
\be D_n(m)=\sum_{k=0}^{r-1}\hat{h}_k\binom{\rho_r-\rho_k}{m}Q_n^{\rho_r-\rho_k}=O(Q_n^{\rho_r-\rho_{r-1}}), \quad m\ge 1,\quad n\to \infty,\la{Ddef}\eu
and
\be
E_n(m)=\binom{\rho_r+1}{m}\left(M(\tilde{\delta}_n;C_0)\right)^\prime Q_n^{\rho_r+1},
\quad m\ge 0,\quad n\to \infty. \la{nush}
\eu
In accordance with \refm[wri],
\be
L_n=\sum_{m\ge 1}(C_n(m)-D_n(m)-E_n(m))Z_n^m-E_n(0).\la{Lnexp}
\eu

We will analyze the factor
$\left(M(\tilde{\delta}_n;C_0)\right)^\prime$ in \refm[nush].
Suppose that \be Z_n=O\left(Q_n^{C_0+\rho_r}\right). \la{Zassump}
\eu Evidently,

$$\left(M(\delta;C_0)\right)^\prime=
\frac{1}{2\pi
i}\int_{-C_0-i\infty}^{-C_0 +i\infty}
\delta^{-s-1}h(s)ds,
$$
where $h(s)$ is
determined by the appropriate choice of (\ref{Mexp1}) - (\ref{Mexp3}).
Under $\tilde{\delta}_n=Q_n(1+Z_n)$, we have
\begin{eqnarray*}
\left(M(\tilde{\delta}_n;C_0)\right)^\prime&=&Q_n^{C_0-1}(1+Z_n)^{C_0-1}
\Big(\int_{-\infty}^{\infty}Q_n^{-i\tau}(1+Z_n)^{-i\tau}h(-C_0+i\tau))id\tau\Big)\\
&=&
Q_n^{C_0-1}(1+Z_n)^{C_0-1}\sum_{l=0}^\infty Z_n^l \Big(\int_{-\infty}^{\infty}Q_n^{-i\tau}{-i\tau\choose l}h(-C_0+i\tau))id\tau\Big).
\end{eqnarray*}

For given $0<C_0<1,\ \rho_r>0$ as in the statement of our Theorem 1,
define
\be m_0=\min\{m\ge 2:(C_0+\rho_r)m> \rho_r+1\},\la{def3}\en
so that $Z_n^{m_0}=o(n^{-1})$.

Next we have
$$E_n(m)=\binom{\rho_r+1}{m} Q_n^{\rho_r+C_0}(1+Z_n)^{C_0-1}\sum_{l=0}^{m_0-2} \eta_l(n)Z_n^l=O(Q_n^{\rho_r+C_0})+ o(n^{-1}),$$$$\quad m\ge 0,\quad n\to \infty,$$
where
\be
\eta_l(n)=\int_{-\infty}^{\infty}Q_n^{-i\tau}{-i\tau\choose l}h(-C_0+i\tau)id\tau=O(1),\quad l\ge 0,\quad n\to \infty
\en
are known constants.
Now the expression for $E_n(m)$ can be written as  $$E_n(m)= Q_n^{\rho_r+C_0}\sum_{l=0}^{m_0-2}\tilde{\gamma}_l(n,m)Z_n^l+ o(n^{-1}),$$
where we denoted
$$
\tilde{\gamma}_l(n,m)=\binom{\rho_r+1}{m}\sum_{p+k=l} \binom{C_0-1}{p}\eta_k(n).
$$

Consequently,
\begin{eqnarray}
\sum_{m\ge 0} E_n(m)Z_n^m& =& \sum_{m= 0}^{m_0-2} E_n(m)Z_n^m  + o(n^{-1})\nonumber\\
&=& Q_n^{\rho_r+C_0}\sum_{k=0}^{m_0-2}Z_n^k\Big(\sum_{m+l=k}\tilde{\gamma}_l(n,m)\Big) +  o(n^{-1})\nonumber\\
&=&
Q_n^{\rho_r+C_0}\sum_{k=0}^{m_0-2}\gamma_k(n)Z_n^k + o(n^{-1}),\la{ENM}
\end{eqnarray}
where
$$
\gamma_k(n)= \sum_{m+l=k}\tilde{\gamma}_l(n,m).
$$

  Set \be W_n(m)=C_n(m)-D_n(m),\ \ m\ge 1 \la{nush1}, \en
Substituting \refm[ENM] and \refm[nush1] into \refm[Lnexp]
produces \be L_n=  \sum_{m=1}^{m_0-1} W_n(m) Z_n^m -
Q_n^{\rho_r+C_0}\sum_{k=0}^{m_0-2}\gamma_k(n)Z_n^k + o(n^{-1}).
\la{juv}\en Our objective will be to show that the condition \be
L_n= o(n^{-1}) \la{tna}\eu is satisfied by $Z_n$ which, in
agreement with the assumption \refm[Zassump], is a polynomial in
$Q_n^{\rho_r+C_0}$ of degree $(m_0-1):$ \be Z_n =
\sum_{k=1}^{m_0-1}\beta_k(n) Q_n^{k(\rho_r+C_0)}, \la{znpol}\eu
with coefficients $\beta_k(n)=O(1)$ that are determined
recursively from \refm[tna]. In fact, under the above form of
$Z_n,$ the main term of $L_n$ becomes a polynomial in
$Q_n^{\rho_r+C_0},$  so that the condition \refm[tna] is satisfied
if \be [Q_n^{k(\rho_r+C_0)}]L_n=0,\  k=1,\ldots,m_0-1,
\la{tna1}\en where $[Q_n^{k(\rho_r+C_0)}]L_n$ denotes the
coefficient of $Q_n^{k(\rho_r+C_0)}$ in $L_n$. From \refm[tna1],
\refm[juv] we derive
$$[Q_n^{(\rho_r+C_0)}]L_n= W_n(1)\beta_1(n)-\gamma_0(n)=0,$$
$$[Q_n^{2(\rho_r+C_0)}]L_n= W_n(1)\beta_2(n) +W_n(2)\beta_1^2(n)-\gamma_1(n)\beta_1(n)=0,\ldots,$$
$$ [Q_n^{k(\rho_r+C_0)}]L_n=[Q_n^{k(\rho_r+C_0)}]\big(\sum_{m=1}^{m_0-1}W_n(m)Z_n^m)\big)- $$$$
[Q_n^{(k-1)(\rho_r+C_0)}]\big(
\sum_{k=0}^{m_0-2}\gamma_k(n)Z_n^k\big)=\beta_k(n)W_n(1)+  R_{k-1}(\beta_1(n), \ldots \beta_{k-1}(n)),$$
$$k=1,\ldots m_0-1,$$ where   $R_{k-1}$ is a polynomial in $\beta_1(n), \ldots \beta_{k-1}(n).$ By induction on $k$, we have $\beta_k(n)=O(1)$ for all
$k\in[1,m_0-1]$.
 Since $R_{k-1}$ does not depend on $\beta_k(n)$
and since $W_n(1)\sim (\rho_r+1) \hat{h}_r\neq 0$ as $n\to\infty,$
our claim is justified. Thus, with $Z_n$ given  by \refm[znpol],
$$L_n=o(n^{-1})+ \sum_{m\ge m_0}W_n(m)Z_n^m\le o(n^{-1})+ O(1)\frac{\vert Z_n^{m_0} \vert}{1-  Z_n}
= o(n^{-1}),$$ by virtue of \refm[def3],\refm[nush1] and because  $Z_n\to 0$.
Consequently,
 we have
$u_n=Z_nQ_n=O(Q_n^{C_0+\rho_r+1})=o(n^{-1}), $
as was desired.

{\bf Step 3} The first two steps, combined with \refm[nmi] produce the following asymptotic formulae for the solutions
$\delta_n^{(i)}:$
\be
\delta_n^{(i)}=\big(\hat{h}_r^{(i)}\big)^{\frac{1}{\rho_r+1}}n^{-\frac{1}{\rho_r+1}}\Big(1+\frac{\sum_{k=0}^{r-1}
\hat{h}_k^{(i)}
\big( \hat{h}_r^{(i)}n^{-1}\big)^{\frac{\rho_r-\rho_k}{\rho_r+1}}+V(z)}{\hat{h}_r^{(i)}}\Big)^{\frac{1}{\rho_r+1}}+o(n^{-1}),
\la{sd}
\en
where $V(z)$ is given by  \refm[jur].
We now observe the following fact stemming from the definition of the set $\Upsilon_r.$
If $\lambda_s,\lambda_t\in \Upsilon_r$ are such that $\lambda_s+\lambda_t\le \rho_r+1,$
then $\lambda_s+\lambda_t\in \Upsilon_r$ as well.
Hence, the binomial expansion of the right hand side of \refm[sd] gives the claimed formulae \refm[sd6].\qed

The following  particular cases of Proposition 1  are of  importance.\\

\noindent {\bf Corollaries} We recall that the set $\tilde{\Upsilon}_r$ is not empty if and only if $\alpha_1=2(\rho_r-\rho_{r-1})\le\rho_r+1.$

$(i)$ {\bf The case when $\rho_r>2\rho_{r-1}$:} Though in this case the set  $\tilde{\Upsilon}_r$ is not necessary empty, \refm[sd] conforms to
 \be
\delta_n^{(i)}=(\hat{h}_r^{(i)})^{\frac{1}{\rho_r+1}}n^{-\frac{1}{\rho_r+1}}+\frac{\sum_{k=0}^{r-1}\hat{h}_k^{(i)}
 (\hat{h}_r^{(i)})^{-\frac{\rho_k}{\rho_r+1}}n^{-\frac{\rho_r-\rho_k+1}{\rho_r+1}}}{\rho_r+1}+o(n^{-1}).
\en
For the proof recall that if $\tilde{\Upsilon}_r\neq
\varnothing$, then $V(z)\sim B_1z^{2(\rho_r-\rho_{r-1})}),$ so that in the expansion  \refm[sd],
$n^{-\frac{1}{\rho_r+1}}V(z)=o(n^{-1})$.

$(ii)$ {\bf The case of a single simple pole:} $r=1, \ \rho_1>0.$
Since $\rho_0=0,$ the condition of Corollary $(i)$ holds and, therefore
\be
\delta_n^{(i)}=(\hat{h}_1^{(i)})^{\frac{1}{\rho_1+1}}n^{-\frac{1}{\rho_1+1}}+\frac{\hat{h}_0^{(i)}
 n^{-1}}{\rho_1+1}+o(n^{-1}),
\en
which recovers the asymptotic formulae for $\delta_n^{(i)}$ in \cite{GSE}. In this regard we note that in
\cite{GSE} a detailed form of the term $o(n^{-1})$ was derived. The latter is not needed for our subsequent study.

$(iii)$ {\bf The case of equidistant simple poles:} $\rho_l=la,\ l=0,\ldots,r, \ a>0$ is a given number.
In this case the condition of Corollary (i) holds only if $r=1.$
For $r\ge 2$, the set
$\Upsilon_r$ consists of  multiples of $a>0$:
$$\Upsilon_r=
\left\{p=aM:\ M\in\BZ_+, M\leq r+\lfloor a^{-1}\rfloor\right\}.$$
 The simple structure of $\Upsilon_r$ allows to write  the expression \refm[sd] for $\delta_n^{(i)}$ in the following form
$$\delta_n^{(i)}=(\hat{h}_r^{(i)})^{\frac{1}{ar+1}}\sum_{l=0}^{r+\lfloor a^{-1}\rfloor}\psi_l^{(i)}z^{al+1}+o(n^{-1}),$$
where $\psi_0^{(i)}=1, \ i=1,2,3.$ We will demonstrate that the  coefficients $\psi_l^{(i)}$ are determined recursively from \refm[condd]. In view of the above representation of $\delta_n^{(i)}$  the condition \refm[condd] can be written as
$$
\hat{h}_r^{(i)}\Big(\sum_{\ell=0}^{r+\lfloor a^{-1}\rfloor}\psi_\ell^{(i)} z^{al}\Big)^{ar+1}-
\sum_{k=0}^{r}\hat{h}_k^{(i)}\Big(\sum_{\ell=0}^{r+\lfloor a^{-1}\rfloor}(\hat{h}_r^{(i)})^{\frac{1}{ar+1}}\psi_\ell^{(i)} z^{al+1}\Big )^{a(r-k)}=o(n^{-1}),
$$
which is equivalent to
$$\hat{h}_r^{(i)}[z^{as}]\Big(\sum_{\ell=0}^{s}\psi_\ell^{(i)}
z^{al}\Big)^{ar+1}-[z^{as}]\sum_{k=0}^{r}\hat{h}_k^{(i)}(\hat{h}_r^{(i)})^{\frac{a(r-k)}{ar+1}}z^{a(r-k)}\Big(\sum_{\ell=0}^{s}\psi_\ell^{(i)}
z^{al}\Big)^{a(r-k)}=0,$$$$s=0,\ldots,r+\lfloor a^{-1}\rfloor .$$
 This gives
\be\label{psirecur} (ar+1)\hat{h}_r^{(i)}\psi_s^{(i)}+q_s(\psi_0^{(i)},\psi_1^{(i)},\ldots,\psi_{s-1}^{(i)})=0,\
s=1,\ldots,r+1,\en where $q_s$ is a polynomial  of
$\psi_1^{(i)},\ldots,\psi_{s-1}^{(i)}$ only. \qed

$(iv)$ {\bf The case where $V(z)$ influences the main term}.\  In all of the previous corollaries
the remainder  $V(z)$ does not influence the  main terms, i.e. terms of order $\ge n^{-1}$,  in the asymptotics of $\delta_n.$
We give now an example, where the impact of $V(z)$ is seen.
Suppose $r=2$ and $2(\rho_2-\rho_1)<\rho_2.$ Then $\alpha_1=2(\rho_2-\rho_1)\in\tilde{\Upsilon}_r,$ and moreover,
 $zV(z)\sim zB_1z^{\alpha_1}\gg n^{-1}, $ as $n\to\infty$.

\nin \subsection{Completion of the proof }

  It is left to use Lemma~\ref{estimates} and Proposition~\ref{deltaexp} to find the asymptotics of the logarithms of each of the three factors in the representation \refm[rep]
  when $\delta=\delta_n^{(i)}$, $i=1,2,3$.

$(i)$
\be
n\delta_n^{(i)}=\big(\hat{h}_r^{(i)}\big)^{\frac{1}{\rho_r+1}}n^{\frac{\rho_r}{\rho_r+1}}+ \sum_{s:\lambda_s\le \rho_r} K^{(i)}_s n^{\frac{\rho_r-\lambda_s}{\rho_r+1}}+ \epsilon^{(i)}_n, \quad i=1,2,3, \la{ndel}
\en
where $\lambda_s\in \Upsilon_r,\ \epsilon_n^{(i)}\to 0.$

$(ii)$ Firstly, by the argument similar to the one for the proof of \refm[arg] we conclude that \be \log f_n^{(i)}(e^{-\delta^{(i)}_n})=\log \nicef^{(i)}(\delta^{(i)}_n)+ \epsilon^{(i)}_n, \quad i=1,2,3. \la{bog}\en
Next, for $l=0,1,\ldots, r$
$$
\big(\delta^{(i)}_n\big)^{-\rho_l}= \big(\hat{h}_r^{(i)}\big)^{\frac{-\rho_l}{\rho_r+1}}n^{\frac{\rho_l}{\rho_r+1}} + \sum_{s:\lambda_s\le \rho_l} K^{(i)}_{s,l} n^{\frac{\rho_l-\lambda_s}{\rho_r+1}}+\epsilon_{n}^{(i)}(l),
$$
where
$\epsilon_{n}^{(i)}(l)=o\left(1\right)$, $i=1,2,3,\ l=1,2,\ldots r$
and where the coefficients $K^{(i)}_{s,l}$ are obtained from the binomial expansion for $\big(\delta^{(i)}_n\big)^{-\rho_l}$, based on \refm[sd6] and the definition \refm[sd2] of the set $\Upsilon_r$.
Consequently, substituting $\delta=\delta^{(i)}_n$ into \refm[prod] gives
$$
\log f_n^{(i)}(e^{-\delta^{(i)}_n})= \sum_{l=0}^{r}\hat{h}_l^{(i)}\big(\hat{h}_r^{(i)}\big)^{\frac{-\rho_l}{\rho_r+1}}n^{\frac{\rho_l}{\rho_r+1}}+
\sum_{l=0}^r \hat{h}_l^{(i)}\sum_{s:\lambda_s\le \rho_l} K^{(i)}_{s,l} n^{\frac{\rho_l-\lambda_s}{\rho_r+1}}+$$
\be \Big(\frac{D(0)}{\rho_r+1}\log n-  \frac{D(0)}{\rho_r+1}\log\hat{h}_r^{(i)}\Big ){\bf 1}(i)+\epsilon_n^{(i)},
\quad i=1,2,3. \la{logfn}
\en

$(iii)$ It turns out that Theorem~4 of \cite{GSE}, which is a local limit theorem, may be adapted to the situation here. We will explain in the proof  below that the reason for it is the fact \refm[gad], which we pointed in our historical remark.
\begin{theorem}({\bf Local limit theorem}).

Let the random variable $U_n$ be defined as in Section 1. Then

\begin{align}
\PP\left(U^{(i)}_n=n\right) &  \
 \sim\frac{1}{\sqrt{2\pi {\rm Var(U_n^{(i)})
 }}}
\sim\frac{1}{\sqrt{2\pi
K_2^{(i)}}}\left(\deli\right)^{1+\rho_r/2}\non\\
& \sim\frac{1}{\sqrt{2\pi
K_2^{(i)}}}\big(\hat{h}_r^{(i)}\big)^{\frac{2+\rho_r}{2(\rho_r+1)}}n^{-\frac{2+\rho_r}{2(\rho_r+1)}},\quad
 n\rightarrow\infty,\quad i=1,2,3,\la{lll}
\end{align}
with constants $K_2^{(i)}$ defined by
$$
K_2^{(1)}= A_r\Gamma(\rho_r+2)\zeta(\rho_r+1),
$$
$$
K_2^{(2)}= A_r(1-2^{-\rho_r})\Gamma(\rho_r+2)\zeta(\rho_r+1)
$$
and
$$
K_2^{(3)}= A_r\Gamma(\rho_r+2).
$$
\end{theorem}

{\bf Proof.} We will sketch the proof that follows the pattern in \cite{GSE}.
Denoting $$\phi_n(\alpha)=\BE \Big(e ^{2\pi i\alpha U_n}\Big),
\quad \alpha\in\BR $$ the characteristic function of
the random variable $U_n$ and setting
$$\alpha_0=(\delta_n^{(i)})^{\frac{\rho_r+2}{2}}\log n,$$ we write
$$
 \PP(U_n=n)=\int_{-1/2}^{1/2}\phi_n(\alpha)
e^{-2\pi in\alpha}d\alpha=I_1+I_2,$$
 where
\be
I_1=\int_{-\alpha_0}^{\alpha_0}\phi_n(\alpha)e^{-2\pi
in\alpha}d\alpha \eu and \be
I_2=\int_{-1/2}^{-\alpha_0}\phi_n(\alpha)e^{-2\pi
in\alpha}d\alpha +\int_{\alpha_0}^{1/2}\phi_n(\alpha)e^{-2\pi
in\alpha}d\alpha. \en
Defining $B_n$ and $T_n$ by \be
  B_n^2=
\Big(\log f_n\big(e^{-\delta}\big)\Big)^{\prime\prime}_{\delta=\delta_n^{(i)}} \ {\rm and}
\quad T_n=-\Big(\log f_n\big(e^{-\delta}\big)\Big)^{\prime\prime\prime}_{\delta=\delta_n^{(i)}}\label{BT}
\end{equation}
for $n$ fixed we have the expansion
 \begin{eqnarray}\label{expand2}
 \non\phi_n(\alpha)e^{-2\pi in\alpha}
&=&
\exp{\left(2\pi i\alpha(\BE U_n-n)-2\pi^2\alpha^2B_n^2+O(\alpha^3)T_n\right)}\nonumber\\
&=&\exp{\left(-2\pi^2\alpha^2B_n^2+O(\alpha^3) T_n\right)},\quad
\alpha\rightarrow0. \label{exp}
 \end{eqnarray}
  By virtue of \refm[prod] and \refm[bog] we derive from \refm[BT] that  the main terms in the asymptotics for $B_n^2$ and   $T_n$ depend on the rightmost pole $\rho_r$ only:
$$B_n^2\sim K_2^{(i)}(\delta_n^{(i)})^{-\rho_r-2}, $$
$$T_n\sim K_3^{(i)}(\delta_n^{(i)})^{-\rho_r-3},$$
where $K_2^{(i)},K_3^{(i)}$ are as in the statement of the theorem.
Therefore, in all three cases, $$B^2_n\alpha_0^2\to \infty, \ T_n\alpha_0^3\to 0, \ n\to \infty.$$
Consequently, $$I_1\sim \frac{1}{\sqrt{2\pi B_n^2}}, \ n\to \infty$$
and it is left to show that $I_2=o(I_1), \ n\to \infty.$
Taking into account that for all three models $\alpha_0=o(\sqrt{\delta_n}),$ we then follow \cite{frgr} splitting  the range $[\alpha_0,1/2]$ of the integral $I_2$ into three subintervals
$[\alpha_0,\delta_n]$,
$[\delta_n,\sqrt{\delta_n}]$, and $[\sqrt{\delta_n},1/2]$.

The proof of Lemma~3 in \cite{GSE} shows that
$$
\phi_n(\alpha)\leq (1+\varepsilon_n)\exp\left(-\frac{V_n(\alpha)}{{\cal M}^{(i)}}\right),\quad \alpha\in \BR,
$$
where $\varepsilon_n\to 0$, ${\cal M}^{(i)}$, $i=1,2,3$, are
positive constants defined in condition $(iii')$ of Theorem~3 in \cite{GSE}, and
$$
V_n(\alpha):= \sum_{k=1}^\infty e^{-k\delta_n}\sin^2(\pi\alpha k).
$$

In the third subinterval, the condition $(III)$  holds, which allows to derive, in the same way as in \cite{GSE}, that the corresponding part of the integral $I_2 $ is $O\left(\delta_n^{\frac{2+\rho_r}{2}+\epsilon}\right)=o(I_1)$.

Regarding the first subinterval, we are able to derive the desired  bound, following the scheme in \cite{frgr}, modified so as to match our setting for the sequence $\{b_k\}$. The subsequent analysis uses an
inequality (3.70) from \cite{frgr}, which we recall in \refm[th] below.
Let $[x]$ and $\{x\}$ denote, respectively,
the integer and fractional parts of a real number $x$, and let
$\Vert x\Vert$ denote the distance from $x$ to the nearest integer,
so
$$
\Vert x\Vert=
\left\{
\begin{array}{l l}
\{x\}&{\rm if \ }\{x\}\leq 1/2;\\
1-\{x\}&{\rm if \ }\{x\}>1/2.
\end{array}
\right.
$$
We then have
\be
\sin^2(\pi x) \geq 4\parallel x\parallel^2. \la{th}
\en

For $\alpha\in [\alpha_0,\delta_n],$ we  get the estimate, as $n\to \infty$:
\begin{eqnarray}
V_n(\alpha)&=&
\sum_{k=1}^n b_ke^{-k\delta_n}\Vert \alpha j\Vert^2\nonumber\\
&\ge& \alpha_0^2 \sum_{k=1}^{1/(2\delta_n)} b_k k^2e^{-k\delta_n}\nonumber\\
&\ge&
 \alpha_0^2C_1\sum_{k=1}^{1/(2\delta_n)}k^2 b_k\nonumber\\
& \sim&
 \alpha_0^2 C_1 \frac{A_{r}}{\rho_r+2}(2\delta_n)^{-\rho_r-2}\la{Ike}\\
&=&
C_2\log^2 n,\nonumber
\end{eqnarray}
$$
C_1=e^{-1/2},\quad C_2=C_1 \frac{A_{r}}{2^{\rho_r+2}(\rho_r+2)}.
$$
In order to show the asymptotic equivalence \refm[Ike], we used the 
Wiener-Ikehara theorem (see \cite{K}, p.122), which for $\{b_k\}$ in our setting reads as follows:
$$\sum_{k=1}^n b_kk^{1-\rho_r}\sim A_rn,\ n\to \infty.$$
By writing $\sum_{k=1}^n b_kk^s=\sum_{k=1}^n (b_k k^{1-\rho_r})(k^{\rho_r+s-1})$,
for $s\geq 0$,
and using summation-by-parts, we obtain
\be\la{Ike2}
\sum_{k=1}^n b_kk^s\sim\frac{A_r}{\rho_r+s}n^{\rho_r+s},\ n\to \infty,
\en
which, with $s=2$, implies \refm[Ike].

For
$\alpha\in [\delta_n,\sqrt{\delta_n}]],$
define, as in \cite{frgr},
\begin{eqnarray*}
{\cal Q}(\alpha)&=&\left\{1\le k\le n: j+1/4\le \alpha k\le j+3/4,
\quad j=0,1,\ldots, \left[
\frac{4\alpha n-3}{4}\right]\right\}\\
&=&\cup_{j=0}^{[\frac{4\alpha n-3}{4}]}{\cal
Q}_{j}(\alpha),
\end{eqnarray*}
where ${\cal Q}_j(\alpha)$ denotes the set of
all integers $k\in [\frac{4j+1}{4\alpha},\frac{4j+3}{4\alpha}].$
\nin Observe that for any $\alpha\in [\delta_n,\sqrt{\delta_n}]$ and
$j\ge 0,$
 the set ${\cal Q}_j(\alpha)$ is not empty, since in this case
 $ \frac{4j+3}{4\alpha}- \frac{4j+1}{4\alpha}\ge 1.$
Now the  aforementioned inequality from \cite{frgr} yields
 \be
V_n(\alpha)\ge \frac{1}{16}\sum_{k\in {\cal Q}(\alpha)}b_ke^{-k\delta_n}=
\frac{1}{16} \sum_{j=0}^{[
\frac{4\alpha n-3}{4}]} \sum_{{k\in {\cal
Q}_j(\alpha)}} b_ke^{-k\delta_n}. \la{99}\en Next, using \refm[Ike2] with $s=0$ and the fact  that $\delta_n\alpha^{-1}\le 1,\ \alpha\in [\delta_n,\sqrt{\delta_n}],$ we estimate the
inner sum in \refm[99], for $j\ge 0:$
\begin{eqnarray*}
\sum_{{k\in {\cal Q}_j(\alpha)}}
b_ke^{-k\delta_n}&=&\sum_{k=[\frac{4j+1}{4\alpha}]}^{[\frac{4j+3}{4\alpha}]}
b_ke^{-k\delta_n}\\
&\ge&\frac{A_{r}}{\rho_r} e^{-[\frac{4j+3}{4\alpha}]\delta_n}\left(\left[\frac{4j+3}{4\alpha}\right]^{\rho_r}-\left[\frac{4j+1}{4\alpha}\right]^{\rho_r}
+ \max\{1,j^{\rho_r}\}o\left(\left(1/\alpha\right)^{\rho_r}\right)\right)\\
&\ge&
\frac{A_{r}}{\rho_r}(4\alpha)^{-\rho_r} e^{-(j+3/4)} \left((4j+3)^{\rho_r}-(4j+1)^{\rho_r}+o(j^{\rho_r})\right)\\
&\ge&
C\delta_n^{-\rho_r/2}e^{-j}\left((4j+3)^{\rho_r}-(4j+1)^{\rho_r}+o(j^{\rho_r})\right),\quad \alpha\in [\delta_n,\sqrt{\delta_n}],
\end{eqnarray*}
where we denoted $C=\frac{A_{r}}{\rho_r}e^{-3/4}4^{-\rho_r}.$
Now continuing \refm[99], we  get with the help of the Euler integral test
\begin{eqnarray*}
V_n(\alpha) &\ge& \frac{1}{16} C\delta_n^{-\rho_r/2}\sum_{j=0}^{[
\frac{4\alpha n-3}{4}]} e^{-j}\Big((4j+3)^{\rho_r}-(4j+1)^{\rho_r}+o(j^{\rho_r})\Big)\\
&\sim&
\frac{1}{16} C\delta_n^{-\rho_r/2}
(e^{3/4}-e^{1/4})\int_0^\infty e^{-x/4}x^{\rho_r}dx\\
&=&C_1 \delta_n^{-\rho_r/2}, \quad C_1>0,
\quad \alpha\in [\delta_n,\sqrt{\delta_n}],\quad n\to \infty.
\end{eqnarray*}
From the preceding analysis it is easily seen that
$e^{-V_n(\alpha)}=o(I_n), \ n\to \infty,$  for all  $\alpha \in [\alpha_0,\sqrt{\delta_n}].$
\qed
\vskip.5cm
Finally, to completely account for the influence of all $r$  poles $\rho_{},\ldots,\rho_1$, we present the sum of the expressions \refm[ndel],\refm[logfn] obtained for the first two  factors $(i)$ and $(ii),$
in the following form:
\begin{eqnarray*}
n\delta_n^{(i)}+ \log f_n(e^{-\delta^{(i)}_n})&=&
\sum_{l=0}^{r}P_l^{(i)}n^{\frac{\rho_l}{\rho_r+1}}+
\sum_{l=0}^r h_l^{(i)}\sum_{s:\lambda_s\le \rho_l} K^{(i)}_{s,l} n^{\frac{\rho_l-\lambda_s}{\rho_r+1}}\\
&&+ \Big(\frac{D(0)}{\rho_r+1}\log n -  \frac{D(0)}{\rho_r+1}\log\hat{h}_r^{(i)}\Big )1(i)+\epsilon_n^{(i)},
\end{eqnarray*}
where $P_l^{(i)}$ denotes the resulting coefficient of $n^{\frac{\rho_l}{\rho_r+1}}.$

{\bf Acknowledgement}

The authors are grateful to a referee who carefully read the paper and made important critical remarks on the proof of Proposition 1.


\begin{thebibliography}{99}

\bibitem{BFHH} Benvenuti, S., Feng B., Hanany, A, He,Y. (2007).
Counting BPS Operators in Gauge Theories -
Quivers, Syzygies and Plethystics. \emph{J.High Energy Physics} 11,050,48 pp.


\bibitem{frgr} G. Freiman, B. Granovsky, (2002). Asymptotic formula for a
partition function of reversible coagulation-fragmentation processes, Isr. J. Math., 130, 259-279.



\bibitem{GSE} Granovsky, B., Stark D. and Erlihson M. (2008).
Meinardus' theorem on weighted partitions:
Extensions and a probabilistic proof.
\emph{Adv. Appl. Math.} \textbf{41} 307-328.


\bibitem{K} Korevaar, J. (2004). \emph{Tauberian Theory}.
Springer.


\bibitem{Kh} Khinchin, A. I., (1960). Mathematical foundations of quantum
statistics, Graylock Press, Albany, N.Y..



\bibitem{LR}
Lucietti, J. and Rangamani, M. (2008).
Asymptotic counting of BPS operators in superconformal field theories. J. Math. Phys. \textbf{49} 30 pp.

\bibitem{M} G. Meinardus, (1954).
Asymptotische Aussagen \"uber Partitionen,
Math. Z. \textbf{59} 388--398.

\bibitem{MW} Madritsch, M. and Wagner, S. (2010).
 A central limit theorem for integer partitions, Monatsh. Math. \textbf{161} 85-114.

\bibitem{TT} Tate, T. (2010). A spectral analogue of the Meinardus theorem on
asymptotics of the number of partitions,
Asymptotic Analysis, \textbf{67}, 1-2, 101-123.




 \bibitem{Pi} Pitman, J. (2006). Combinatorial stochastic processes. \emph{Lecture Notes in Mathematics, 1875}.

 \bibitem{V1} Vershik, A. (1996). Statistical mechanics of combinatorial partitions and their limit
  configurations. \emph{Funct. Anal. Appl.} \textbf{30}, 90-105.



\end{thebibliography}
\end{document}